\documentclass[11pt,a4paper,psamsfonts]{amsart}
\usepackage{amssymb,amscd,amsmath,amsthm}

\chardef\bslash=`\\ 

 \theoremstyle{plain} 
\newtheorem{theorem}{Theorem}[section]
\newtheorem{corollary}[theorem]{Corollary}
\newtheorem{lemma}[theorem]{Lemma}
\newtheorem{proposition}[theorem]{Proposition}

\theoremstyle{remark}
\newtheorem{remark}[theorem]{Remark}
\newtheorem{example}[theorem]{Example}
\theoremstyle{definition}
\newtheorem{definition}[theorem]{Definition}

\newcommand{\thmref}[1]{Theorem~\ref{#1}}

\newcommand{\proref}[1]{Proposition~\ref{#1}}

\def\K{\mathcal K}

\def\inv{^{-1}}

\newcommand{\matr}[2]{\left(\begin{array}{cc} 1&#1 \\ 0&#2
  \end{array}\right)}
\newcommand{\smatr}[2]
 {\left({\mbox{\tiny{$\begin{array}{cc} 1&#1 \\ 0&#2\end{array}$}}}\right)}
\numberwithin{equation}{section}
\def\adel{\mathcal A_f}

\def\basl{\backslash}
\def\Go{\Gamma_0}
\def\Ga{\Gamma}
\def\hekn{\mathcal H (N,\Gamma_0)}
\def\hekg{\mathcal H (\Gamma,\Gamma_0)}
\def\rels{\mathfrak h}
\begin{document}

\title{Hecke algebras of semidirect products}
\author[M. Laca]{Marcelo Laca*}
\address{Department of Mathematics, University of M\"{u}nster, 48149 M\"unster, Germany.}
\email{laca@math.uni-muenster.de}
\thanks{* ) Supported by the Deutsche Forschungsgemeinschaft [SFB 478]}
\author[N. S. Larsen]{Nadia S. Larsen**}
\address{Department of Mathematics, University of Copenhagen, 
Universitetsparken 5, DK-2100 Copenhagen \O, Denmark.}
\email{nadia@math.ku.dk}
\thanks{**)Supported by the Danish Natural Science Research Council.}
\subjclass{46L55}

\begin{abstract} 
We consider group-subgroup pairs in which the group is a semidirect product
and the subgroup is contained in the normal part. We give conditions for the pair 
to be a Hecke pair and we show that the enveloping Hecke algebra and Hecke
$C^*$-algebra are canonically isomorphic to semigroup crossed products,
generalizing earlier results of Arledge, Laca and Raeburn and of Brenken.
\end{abstract}
\date{June 25, 2001}

\maketitle

\section*{Introduction}
A group--subgroup pair $(\Ga, \Go)$ is called a {\em Hecke pair} if 
 each double coset $\Go\gamma \Go$ is
a finite union of left cosets $\gamma_i \Gamma_0$. By considering inverses, one
sees that it is equivalent to require that
each double coset is the union of finitely many right cosets.
This situation is variously referred to by saying that 
$\Go$ is a {\em Hecke subgroup} or an {\em almost normal} subgroup of $\Ga$. See
\cite{kri} for more details.

To each Hecke pair one associates a {\em Hecke *-algebra} $\mathcal H(\Gamma, \Gamma_0)$, by endowing 
the set of $\Gamma_0$-biinvariant complex valued
functions on $\Gamma$ that are supported in finitely many double cosets with a convolution product 
 defined by the (finite) sum
\[
f\ast g(\gamma)=\sum_{\gamma_1\in \Gamma_0 \basl \Gamma} f(\gamma
\gamma_1^{-1}) g(\gamma_1),
\]
and an involution defined by
$f^*(\gamma)=\overline{f(\gamma^{-1})}$. The characteristic
function of $\Gamma_0$ is the identity element.
Hecke algebras have been around for quite some time but their C*-algebraic version
has only recently begun to get some attention, initially in the work of Binder \cite{bin} and Bost and Connes
\cite{bos-con}, and more recently in that of Hall \cite{hal} and Tzanev \cite{tza}.
One defines the {\em reduced Hecke C*-algebra} $C^*_r(\Ga,\Go)$ to be the closure of $\hekg$ in the natural
representation on $\ell^2(\Go\basl \Ga)$ obtained using
the above convolution formula with $g \in \ell^2(\Go\basl \Ga)$. 
A universal version $ C^*_u(\Ga,\Go)$, through which all representations of $\hekg$ on Hilbert space factor, is
not so readily obtained because in some cases there are generators that do not have uniformly bounded norms in all
*-representations, see \cite[Example 2.3.3]{hal}. To some extent this difficulty can be sidestepped by restricting
the attention to representations that do not increase a certain Banach algebra norm on $\hekg$, see \cite{tza}.

The present work is motivated by the realization of various
Hecke C*-algebras as semigroup crossed products. Following \cite{bcalg,alr}, but especially 
guided by \cite{bre}, we are first led to consider 
the semigroup of elements of $\Ga$ whose left coset is right invariant --and hence a double coset. 
We observe that if this semigroup generates all of $\Ga$ by left quotients, then the problem with unbounded
norms does not arise because $\hekg$ has a generating set consisting of isometries,
providing a simple proof of a result of Hall concerning the existence of 
the universal Hecke C*-algebra \cite[Corollary 4.6]{hal}.

Our main results are about pairs consisting of a semidirect product group $\Ga = N\rtimes G$ 
and a subgroup $\Gamma_0\subset N$.
We give sufficient conditions on the action of $G$ on $N$ for $(\Ga,\Go)$ to be a Hecke pair, and we
realize the Hecke algebra, and C*-algebra, as semigroup crossed products. These results unify and
generalize previous realizations of Hecke C*-algebras as semigroup crossed products obtained in 
\cite{bcalg,alr,bre,LarsenR_MathScand}. 

Our approach relies freely on the basic results from Section 1 of
\cite{bre} where the generality is appropriate for our purposes, but the techniques used to obtain
our main results are independent of the rest of \cite{bre}.
 Indeed, we are able to give simpler proofs of more general
statements. Notably, by considering a different set of generating monomials, inspired in those of \cite{alr},
we have done away with the quasi-lattice structure on $G$.
Perhaps more significantly, we have removed the assumption that $\Go$
be normal in $N$; instead we assume that $(N,\Go)$ is a Hecke pair, and
replace the group algebra $\mathbb C(N/\Go)$ with the Hecke algebra $\hekn$.

We would like to thank Iain Raeburn and George Willis for several helpful 
conversations.

\section{Hecke algebras of semidirect products}
\begin{definition} (cf. \cite{bos-con})
For each $\gamma \in \Gamma$, let $R(\gamma)$ (respectively, $L(\gamma) $) be the cardinality of the
image of $\Gamma_0\gamma\Gamma_0$ in $\Gamma_0 \basl \Gamma$ (respectively, in $\Ga / \Go$). It is clear
that $L(\gamma) = R(\gamma\inv)$. 
\end{definition}
\noindent{\bf Notation.} Throughout we denote by $[\gamma]$ the characteristic function 
of the double coset of an element $\gamma \in \Gamma$.

\begin{remark}\label{prod-formula}
The Hecke algebra $\hekg$ is the linear span of $\{ [\gamma] : \gamma \in \Ga\}$. 
The product of two elements of this spanning set is given by
$([\gamma]  [\delta]) (x) = |\Go \basl (\Go \gamma \inv \Go x\cap \Go \delta \Go)|$,
from which one easily sees that the support of $[\gamma]  [\delta]$ is $\Go \gamma \Go \delta \Go$.
\end{remark}
{}From theorems 1.2 and 1.4 of \cite{bre} we know that the set 
\[
\Sigma := \{\sigma \in \Gamma : L(\sigma) = 1\}
\]
is a semigroup on which $R$ is multiplicative,  and that the map
\[
 \sigma \in \Sigma  \ \mapsto \  W_\sigma  
:= R(\sigma)^{-1/2} [\gamma] \in \mathcal H(\Gamma,\Gamma_0) 
\]
is a representation of $\Sigma$ by isometries in $\mathcal H(\Gamma, \Gamma_0)$.
We will need the following supplement of this. 

\begin{lemma}\label{supplement}
For every $\gamma\in \Gamma$ and $\sigma \in \Sigma$, 
\begin{enumerate}
\item[(i)]the product $[\gamma][\sigma]$ is a positive integer
multiple of $[\gamma \sigma]$, and $L(\gamma \sigma) = L(\gamma)$; 
\item[(ii)] the product $[\sigma]^* [\gamma]$ is a
positive integer multiple of
$[\sigma \inv\gamma]$, and $ R(\sigma\inv \gamma) = R(\gamma)$.
\end{enumerate}
\end{lemma}
\begin{proof}
Clearly $\sigma \in \Sigma $ if and only if $\Gamma_0 \sigma  \Gamma_0 = \sigma \Gamma_0 $, 
from which it follows that
$\Gamma_0 \gamma \sigma \Gamma_0 = \Gamma_0 \gamma \Gamma_0 \sigma \Gamma_0 = 
(\sqcup_{i=1}^{L(\gamma)} \gamma_i\Gamma_0) \sigma\Gamma_0 =
\sqcup_{i=1}^{L(\gamma)} \gamma_i \sigma \Gamma_0$.
Since the support of $[\gamma][\sigma]$ is $\Gamma_0
\gamma \Gamma_0 \sigma \Gamma_0$, we obtain (i), from which (ii) follows
on  considering adjoints and replacing $\gamma$ by $\gamma\inv$. 
\end{proof}

\begin{proposition}\label{prop1}
Suppose $(\Ga , \Go)$ is a Hecke pair such that $\Sigma\inv \Sigma = \Gamma$.
Then $\{W_\tau^* W_\sigma^{}: \tau, \sigma\in \Sigma\}$ is a linear basis
of $\hekg$, and the universal C*-enveloping algebra $C^*_u(\Ga, \Go)$
of $\hekg$ is a quotient of the universal C*-algebra of isometric
representations of $\Sigma$.
\end{proposition}
\begin{remark} Since the semigroup $\Sigma$ is cancellative, because it is contained in 
$\Ga$, by assuming that $\Sigma\inv \Sigma = \Gamma$ (or equivalently that
$\Sigma \gamma  \cap \Sigma \gamma' \neq \emptyset$ for every pair
$\gamma,\gamma' \in \Ga$) we are saying that $\Sigma$ is an 
Ore semigroup with $\Ga$ its group of left quotients.
The proposition asserts in particular that $C^*_u(\Ga, \Go)$ exists, proved in \cite[Corollary
4.6]{hal}, see also \cite[Proposition 2.8]{bre}.
\end{remark}

\begin{proof}
Let $\tau, \sigma \in \Sigma$; from the lemma we know that
 \[
[\tau ]^*  [\sigma]
=
K_{ \sigma, \tau}[\tau\inv \sigma ]
\]
where $K_{\sigma,\tau}$ is a positive integer. 
Hence, up to positive constants, the products $W_\tau^* W_\sigma^{}$ are the characteristic  
functions of all the double cosets, so their linear span is $\mathcal H(\Gamma,\Gamma_0)$.
{}From this it already follows that $\mathcal H(\Gamma, \Gamma_0)$ is generated by isometries
and that the universal C*-enveloping algebra
for $\mathcal H(\Gamma, \Gamma_0)$ exists and is a quotient of the universal C*-algebra of
isometric representations of $\Sigma$. The key here is of course that isometries have norm $1$ in any
nontrivial Hilbert space representation, so the basic conditions for the existence of the universal C*-algebra
(see e.g. \cite{bla}) are satisfied.

In order to conclude that the spanning set is a linear basis of $\mathcal H(\Gamma, \Gamma_0)$
we have to prove that
$ W_\beta^* W_\alpha^{}= W_\tau^* W_\sigma^{}$ whenever $\beta\inv \alpha = \tau\inv \sigma$, 
with $\alpha, \beta, \sigma$ and $\tau$ in $\Sigma$. 
For this it suffices to show that the expression
\begin{equation}\label{koverR}
\frac{K_{\sigma,\tau} }{R(\tau)^{1/2} R(\sigma)^{1/2}}
\end{equation}
depends only on the left quotient $\tau\inv \sigma$. 
 Since $\Sigma$ is directed, that is, $\Sigma \gamma \cap \Sigma \gamma ' \neq \emptyset$,
it is enough to show that for each $\gamma \in \Sigma$
the value of the above expression does not change when we replace $\sigma$ and $\tau$ by
$\gamma \sigma$ and $\gamma \tau$, respectively, for which it suffices to show that
$K_{\gamma \sigma,\gamma \tau} = R(\gamma) K_{\sigma,\tau}$ because $R$ is  multiplicative  on
$\Sigma$. In order to prove this we first compute $K_{\sigma,\tau}$ by evaluating the
convolution product at $ \tau\inv \sigma$ using Remark \ref{prod-formula}:
\begin{align*}
K_{\sigma,\tau} & = ([\tau\inv]  [\sigma]) (\tau\inv \sigma)\\
 & =  \vert\Gamma_0 \backslash (\Gamma_0\tau\Gamma_0 \tau\inv \sigma \cap
\Gamma_0\sigma\Gamma_0)\vert  \\ & =  \vert\Gamma_0 \backslash (\tau\Gamma_0 \tau\inv \sigma \cap
\sigma\Gamma_0)\vert \\ & =  \vert \Gamma_0 \backslash (\tau\Gamma_0 \tau\inv \cap
\sigma\Gamma_0 \sigma\inv)\vert.
 \end{align*}
The claim then follows on applying the basic fact 
$\vert A \backslash C\vert = \vert A \backslash B\vert \, \vert B \backslash C\vert$, valid
for every group inclusion $A \subset B \subset C$,
to the inclusions
\[ 
\Gamma_0 \subset \gamma \Gamma_0 \gamma\inv 
\subset \gamma\tau\Gamma_0 \tau\inv\gamma\inv \cap \gamma\sigma\Gamma_0 \sigma\inv\gamma\inv,
\]
and observing that
$R(\gamma) = \vert\Gamma_0 \backslash (\gamma \Gamma_0 \gamma\inv )\vert $ and that
\[
\vert (\gamma \Gamma_0 \gamma\inv )\backslash 
(\gamma\tau\Gamma_0 \tau\inv\gamma\inv  \cap \gamma\sigma\Gamma_0 \sigma\inv\gamma\inv)\vert
= 
\vert \Gamma_0 \backslash 
(\tau\Gamma_0 \tau\inv \cap \sigma\Gamma_0 \sigma\inv)\vert.
\]
\end{proof}

We pause to recall that when $\alpha$ is an action of 
a semigroup $S$ by endomorphisms of a unital C*-algebra $A$, one defines a {\em covariant representation} to be a
pair $(\pi,V)$ in which $\pi$ is a unital representation of $A$ and $V$ is a representation of $S$ by isometries
 such that $\pi (\alpha_s(a)) = V_s \pi(a) V_s^*$ for every $a\in A$ and
$s\in S$. The associated {\em semigroup crossed product} is the C*-algebra generated by a universal covariant representation,
which, by definition, is one through which every covariant representation factors, see \cite{quasilat, mur} for more details. 
The analogous concept in the category of *-algebras is defined in the obvious way and we refer to it as 
the *-algebraic semigroup crossed product. 

\begin{remark}\label{trick}
Since $W$ is a semigroup of isometries, there is a semigroup of endomorphisms 
(of $\hekg$ and of $C^*_u(\Ga,\Go)$) defined by $\alpha_s(b) := W_s b W_s^*$.
Let $B$ denote the *-subalgebra of $\hekg$ generated by the projections 
$W_sW_s^*$ with $s\in \Sigma$, in other words, $B$ is the linear span of the products of such projections. 
Since $W_\sigma (\prod_j W_{s_j}W_{s_j}^*) W_\sigma^* =\prod_j W_{\sigma s_j} W_{\sigma s_j}^*$, the endomorphisms 
leave $B$ invariant, indeed, $B$ is the smallest invariant unital *-subalgebra.   
It is easy to show that the associated semigroup crossed product $B\rtimes
\Sigma$ maps onto the Hecke algebra in a canonical way. However, there is in principle no reason why this map should be
injective in general, because the covariance relation (which takes place in the algebra $B$) might not capture all the
relations satisfied by the generators of the Hecke algebra. 
By the same token, if $N$ is a normal subgroup of $\Ga$ and $\Go$ is normal in $N$, there is
a surjective homomorphism of $\mathbb C (N/\Go) \rtimes \Sigma$ onto $\hekg$. This shows that the conclusion of
\cite[Proposition 2.10]{bre} holds without any assumptions on least upper bounds.
\end{remark}

When $\Ga$ is a semidirect product we will be able to say more.
 Recall that if the group $G$ acts on the group 
$N$ by automorphisms $\psi_g$, one defines the semidirect product group $\Ga = N\rtimes G$ by endowing the
set $N\times G$ with the twisted product $(n,g) (m,h) = (n\psi_g(m), gh)$. One can then view $G$
and $N$ as subgroups of $\Ga$, with $N$ normal, and write the automorphisms $\psi_g$ as inner automorphisms
obtained by conjugation by elements of $G$, namely, $\psi_g(m) = g m g\inv$.

\begin{proposition} \label{mainprop}
Let $\Gamma := N\rtimes G$ be a semidirect product group and let $\Go$ be a subgroup of $N$.
Suppose that $S$ is a subsemigroup of $G$ such that $S\inv S = G$
and $s\inv \Go s \subset \Go$ for every $s\in S$. 
Then $(\Ga,\Go)$ is a Hecke pair if and only if $(N,\Go)$ is
a Hecke pair  and $\vert s\inv \Gamma_0 s \basl \Gamma_0 \vert < \infty$ for every $s \in S$. Furthermore,
 $R(s) = \vert s\inv \Gamma_0 s \basl \Gamma_0 \vert$ for every $s\in S$, and the identity map
gives an embedding $\hekn \hookrightarrow \hekg$.
\end{proposition}
\begin{proof}
  The assumption $s\inv \Go s \subset \Go$ 
says that the action of $S\inv$ leaves $\Go$ invariant. This is also equivalent to saying that 
 $S\subset \Sigma$, because $s\inv \Go s \subset \Go$ implies $\Go s\Go = s\Go$.

Assume first that $\Go $ is a Hecke subgroup of $\Ga$.  Then it is clearly also a 
Hecke subgroup of $N$, and for each $s\in S$ one has $\vert s\inv \Gamma_0 s \basl \Gamma_0 \vert =
\vert \Go \basl s \Go s\inv \vert = \vert \Go \basl \Go s \Go s\inv \vert =
\vert \Go \basl (\Go s \Go )\vert = R(s) < \infty$.

Assume now that $\Go $ is a Hecke subgroup of $N$ and that
 $\vert s\inv \Gamma_0 s \basl \Gamma_0 \vert < \infty$ for every $s\in S$.
Let $\gamma\in \Ga$ and write $\gamma = t\inv n s$ with $t,s\in S$ and $n\in N$,
(an easy argument shows this is always possible because $S\inv S = G$). In order to show that 
$\vert\Go \basl \Go t\inv n s \Go\vert$ is finite, we write:
\begin{align*}
\Go t\inv n s \Go & = \Go t\inv \Go n \Go s \Go\\
& = \cup_j \Go t\inv \Go n \Go \gamma_j  \qquad\quad 
\text{because } s \Go = \Go s \Go =  \sqcup_{j=1}^{R(s)} \Go \gamma_j\\
& = \cup_j \cup_i  \Go t\inv \Go n_i \gamma_j  \qquad 
\text{ because }\Go n \Go = \sqcup_{i=1}^{R(n)} \Go n_i\\
& = \cup_j \cup_i  \Go t\inv n_i \gamma_j ,
\end{align*}
which is a union of at most $R(s) R(n)$ right cosets. This proves that $(\Ga,\Go)$ is a Hecke pair,
and hence that $R(s) = \vert s\inv \Go s \basl \Go \vert$, as shown above. 
Finally, the operations defined on $\hekg$ restrict to those on $\hekn$,
so the identity map gives a
 unital embedding of *-algebras $\hekn \hookrightarrow \hekg$.
\end{proof}
\begin{remark} 
The set $S_0 :=\{ s\in G: [\Gamma_0 : s\inv \Gamma_0 s]<\infty\}$ is itself always a semigroup, 
in fact, the largest possible subsemigroup of $G$ that can satisfy the conditions of hypothesis. 
If $S_0\inv S_0 = G$, the proposition certainly applies with $S = S_0$, 
but in some situations one might be better off with a proper subsemigroup $S$ of $S_0$ that
happens to generate $G$ by left quotients too. 
\end{remark}

\begin{theorem} \label{mainthm}
Let $\Gamma := N\rtimes G$ be a semidirect product group and let $\Go$ be a subgroup of $N$. Suppose that $S$ 
is a subsemigroup of $G$ such that $S\inv S = G$ and $\vert s\inv \Go s \basl \Go \vert < \infty$ for every $s\in S$.
If  $(N,\Go)$ is a Hecke pair, then $(\Ga,\Go)$ is a Hecke pair, with $R(s) = \vert s\inv \Go s \basl \Go \vert$, and

\begin{itemize}
\item[(i)] 
the elements of $\hekg$ given by  $\mu_s := R(s)^{-1/2} [s]$ for $s\in S$ and  $e(x):= [x]$ for $x\in N$ satisfy the relations
\begin{enumerate}
\smallskip
\item[($\rels_1$)]\   $\mu_s^*\mu_s= I$ and $\mu_s\mu_t=\mu_{st}$,  for $s, t\in S$;
\smallskip
\item[($\rels_2$)]\  $e(1)=I$, $e(x)^*= e(x\inv)$, and $e(x) \, e(y) = [x]  [y] $ (as in $\mathcal H (N,\Gamma_0)$);
and
\smallskip
\item[($\rels_3$)]\  $\mu_s e(x)\mu_s^*= R(s)\inv {\sum _{\{ 
[y] \,:\,  [s\inv y s] = [x] \}}} e(y)$, for $s\in S$ and $x \in  N$;
\smallskip
\end{enumerate}
and the set $\{\mu_t^* e(x) \mu_s : \, s,t \in S \text{ and } x\in N\}$
is a linear basis of $\hekg$;

\item[(ii)] there is an action $\alpha$ of $S$  by injective corner endomorphisms of the Hecke *-algebra 
$\mathcal H(N,\Gamma_0)$ defined by
\[
\alpha_s(e(x)) := R(s)\inv {\sum_{\{ 
[y] \,:\,  [s\inv y s] = [x] \}}} e(y),
\] 
 and the Hecke *-algebra $\hekg$ is canonically isomorphic
to the *-algebraic semigroup crossed product $\hekn
\rtimes_\alpha S$, and to the universal unital *-algebra with presentation ($\rels_1,\rels_2, \rels_3$).
\end{itemize}
\end{theorem}
\begin{proof}
The first assertion is from \proref{mainprop}.
A somewhat tedious but straightforward computation of convolution products in $\hekg$ shows that the 
$\mu_s $ and the  $e(x)$ satisfy the relations ($\rels_1$, $\rels_2$, $\rels_3$). 
Notice that relation ($\rels_1$) simply says that $\mu$ is an isometric representation of $S$
(see Theorem 1.4  of \cite{bre}) and
relation ($\rels_2$) is just a restatement of the embedding of $\hekn$ in $\hekg$.
In order to see that $\{\mu_t^* e(x) \mu_s : \, s,t \in S \text{ and } x\in N\}$
is a linear basis of $\hekg$, notice first that the support of  $\mu_t^* e(x) \mu_s$
is the set $\Go t\inv \Go x \Go s \Go = \Go t\inv x s \Go$, and that since $\Ga = S\inv N S$,
 all double cosets arise
as supports. Thus the products $\mu_t^* e(x) \mu_s$ linearly span $\hekg$, and 
we only have to prove that they are linearly independent.  Different supports
are disjoint, so to prove linear independence it suffices to show that if two such products have the same
support then they coincide, which we do in the following lemma.
\begin{lemma}
Let $t, s, \tau, \sigma \in S$ and $ m,n \in N$, and assume that 
$[\tau\inv m \sigma] = [t\inv n s]$. Then
$\mu_\tau^* [m] \mu_\sigma^{} = \mu_t^* [n] \mu_s^{}$.
\end{lemma}
\begin{proof}
First take quotients modulo $N$ in $[\tau\inv m \sigma] = [t\inv n s]$ to see that $\tau\inv \sigma = t\inv
s$.  Since $S$ is directed there exist $r$ and $\gamma$ in $S$ 
such that $rt = \gamma \tau$ and $rs = \gamma \sigma$. 
Since $S\subset \Sigma$, Lemma \ref{supplement} implies that the function $R$ is invariant under multiplication 
by elements of $S\inv$ on the left, so
\[
R(\sigma) = R(\tau \inv \sigma ) = R(t\inv s) = R(s).
\]
Since  $R$ is multiplicative on $S$ and $R(\gamma \sigma) = R(rs)$, we also have $R(r) = R(\gamma)$.
{}From the assumption $\Go \tau\inv m \sigma \Go = \Go t\inv n s \Go$ it follows at once that
\begin{align}\label{assumption}
\Go \tau\inv \Go m \Go \sigma \Go = \Go t\inv \Go n \Go s \Go, 
\end{align}
at which point it becomes convenient to let $x := rt = \gamma \tau$ and $y:= rs = \gamma \sigma$, and to write
 $\tau, \sigma, t$ and $ s$ in terms of $x, y, r$ and $\gamma$.
Multiplication of \eqref{assumption}  on the left by $x$ and on the right by $y\inv$ yields
\[
(x \Go x\inv ) (\gamma \Go m \Go \gamma\inv) (y \Go y\inv)
=
(x \Go x\inv ) (r \Go n \Go r\inv) (y \Go y\inv).
\]
Using the multiplication rules in $\hekg$ we see that
\begin{itemize}
\item[] $(x \Go x\inv )$ is the support of $\mu_x \mu_x^*$,
\item[] $(y \Go y\inv )$ is the support of $\mu_y \mu_y^*$,
\item[] $(\gamma \Go m \Go \gamma\inv)$ is the support of $\mu_\gamma [m] \mu_\gamma^*$,
and 
\item[] $(r \Go n \Go r\inv)$ is the support of $\mu_r [n] \mu_r^*$.
\end{itemize}
It follows that
\[
\mu_x \mu_x^* \mu_\gamma [m] \mu_\gamma^* \mu_y^{}\mu_y^* = \mu_x \mu_x^* \mu_r[n] \mu_r^* \mu_y^{}\mu_y^*,
\]
because the supports match and since $R(\gamma) = R(r)$ the coefficients also coincide.
Now simply  multiply  by $\mu_x^*$ on the left,  by $\mu_y$ on the right, and simplify
$\mu_x^*\mu_x\mu_x^* \mu_\gamma = \mu_\tau^*$,  and other similar products, to see that
$\mu_\tau^* [m] \mu_\sigma^{} = \mu_t^* [n] \mu_s^{}$. This finishes the proof of the lemma
and thus of part (i).
\end{proof}

It is possible to give a direct computational proof of part (ii)
along the lines of \cite{bcalg,alr,bre}, but it is more efficient to
use the Hecke algebra as in Remark \ref{trick}. 
To this effect notice that, because of relation ($\rels_3$) with $x\in \Go$,  the projection $\mu_s \mu_s^*$ is
contained in $\hekn$ for every $s\in S$. Since $\mu_s$ is a semigroup of isometries by relation ($\rels_1$), the left
hand side of ($\rels_3$) defines a semigroup of injective corner endomorphisms: the image of $\alpha_s$ is the corner
$\mu_s\mu_s^* \hekn \mu_s\mu_s^ *$, on which $\alpha_s$ has an 
inverse defined by $X \mapsto \mu_s^* X \mu_s$.

A standard reinterpretation of 
the relations ($\rels_1$, $\rels_2$, $\rels_3$)  shows that they are in fact a presentation of 
the *-algebraic semigroup crossed product $\hekn \rtimes_\alpha S$.
Since the $\mu_s$'s and the $e(x)$'s in the Hecke algebra satisfy the relations, 
there is a canonical homomorphism of $\hekn \rtimes_\alpha S$ to $\hekg$, which
we must show is an isomorphism. In order to avoid confusion we temporarily 
denote the generators of the universal *-algebra with the given presentation
by $\tilde\mu(s)$ and $\tilde e(x)$ while maintaining $\mu_s$ and
$e(s)$ for the generators of the Hecke algebra.  The linear span of products of the form 
$\tilde\mu_t^* \tilde e(x) \tilde\mu_s $ is multiplicatively closed, and so must be all of  $\hekn
\rtimes_\alpha S$, see for example Remark 1.6 of \cite{minautex}. But the image of this spanning set is the
linear basis $\{\mu_t^* e(x) \mu_s : s,t\in S, \ x\in N\}$ of the Hecke algebra obtained in part (i).
Since a spanning set is mapped to a linear basis, the former is also a
linear basis and the map is an isomorphism. This finishes the proof of (ii)
and of the theorem.
\end{proof}

An immediate application of \thmref{mainthm} and the corresponding universal properties gives the following
result about C*-algebras.
\begin{theorem} \label{maincstarthm}
Under the  assumptions of the preceding theorem,
the universal Hecke C*-algebra $C^*_u(\Ga,\Go)$ exists 
if and only if $C^*_u(N,\Go)$ exists,  in which case the isomorphisms of part (ii) 
 extend to C*-algebra isomorphisms:
\[C^*_u(\Ga,\Go) \cong C^*_u(N,\Go) \rtimes_\alpha S \cong C^*_{u} (\rels_1, \rels_2, \rels_3).
\]
\end{theorem} 
Notice that, by \cite[Theorem 2.1.1.]{minautex}, 
 the first isomorphism implies that  $C^*_u(N,\Go)$ embeds in $C^*_u(\Ga,\Go)$.

\begin{corollary}\label{normal}
If $\Go$ is normal in $N$, then relation {\em ($\rels_2$)} defines the group algebra of $N/\Go$ so one has
\[
\hekg \cong \mathbb C (N/\Go) \rtimes_\alpha S
\quad \text{ and} \quad
C^*_u(\Ga,\Go) \cong C^*(N/\Go) \rtimes_\alpha S.
\]
\end{corollary}

\begin{remark}
By the dilation/extension theorems of \cite{minautex} the representation theory of 
$\hekg$ is equivalent to that of a C*-dynamical system  $(A,G,\beta)$ in which $A$ is the
inductive limit $\lim_S (\hekn, \alpha_s)$ and $\beta $ is an action of $G$ by automorphisms of $A$
(we refer to the appendix of \cite{lar-rae} for the categorical statement of this equivalence). This is
particularly  useful when $(N,\Go)$ is a Gelfand pair, that is to say, when $\hekn$ is commutative,
 since it then allows one to apply results about  transformation group C*-algebras to Hecke C*-algebras. 
See for instance \cite{primbc} and \cite{lar-rae}.
\end{remark}

\section{Examples}
\begin{example} \thmref{mainthm} applies to the Bost-Connes Hecke C*-algebra from \cite{bos-con} and to
the Hecke C*-algebras of the inclusions
\[
\matr{\mathcal O}{1} \subset \matr{\mathcal K}{\mathcal K^*}
\]
where $\mathcal K$ is a number field with ring of integers $\mathcal O$.
These have already been shown to be semigroup crossed products in \cite{bcalg} and \cite{alr}, and
also in \cite{bre} under the additional assumption that $\mathcal K$ is of class number $1$.
For the case when  $\mathcal O$ is one of the principal subrings of a global field $\mathcal K$
 see \cite{LarsenR_MathScand}.
\end{example}

\begin{example} We briefly recall now the situation discussed in Example 4.3 of \cite{bre}.
Let $G= GL^+_d(\mathbb Q)$ be the group of rational $d\times d$ matrices with
positive determinant, and let $S = G \cap M_d(\mathbb Z)$ be the subsemigroup of matrices having integer entries.
Then $G$ acts on $N = \mathbb Q^d$ by $\psi_g(n) =(g^t )\inv n $, where $g^t$ indicates the transpose of
$g$ (we choose this action instead of the conventional action used in \cite{bre} in order to
avoid having to switch from $S$ to $S\inv$ later on).
Let $\Gamma := N\rtimes G = \mathbb Q^d\rtimes GL^+_d(\mathbb Q)$, and consider the subgroup $\Gamma_0 =
\mathbb Z^d \subset \mathbb Q^d = N$.  It is clear that if $s\in S$ then
$\psi_{s\inv} (\Gamma_0) = s^t\mathbb Z^d  \subset \mathbb Z^d$ and since
$| \Gamma_0 / \psi_{s\inv} (\Gamma_0)| = \det (s^t) < \infty$, we have a finite order inclusion 
for every $s\in S$. 

{}From \thmref{mainthm} we deduce that $(\mathbb Q^d\rtimes GL^+_d(\mathbb Q), \,\mathbb Z^d )$
is a Hecke pair and that its Hecke C*-algebra is isomorphic to the semigroup crossed product 
\[
C^*_u(\mathbb Q^d\rtimes GL^+_d(\mathbb Q),\mathbb Z^d) \cong C^*(\mathbb Q^d /\mathbb Z^d) \rtimes_\alpha S.
\]
\begin{remark} 
Theorem 2.11 of \cite{bre} suffices for the existence of this universal Hecke
C*-algebra, but it does not give the semigroup crossed product structure, 
because in this case the normalizer of $\Gamma_0$ is not normal in $\Gamma$.
\end{remark}

Next we shall use a multivariable version of the arguments of Section 3 of \cite{minautex} to obtain
an explicit realization of the minimal dilated system whose crossed product contains the Hecke C*-algebra 
as a full corner. Denote by $\adel$ the ring of finite adeles over 
$\mathbb Q$, with $\mathcal R = \prod_p \mathbb Z_p$ the 
maximal compact subring. Corresponding to the natural action $(g,x) \mapsto gx$ of $G$ on $\adel^d$ there is
an action $\beta$ of $G$ on $C_0(\adel^d)$  given by $ (\beta_g F)(x) = F(g\inv x)$ for $F\in C_0(\adel^d)$. 
 
\begin{proposition} The Hecke C*-algebra
$ C^*_u(\mathbb Q^d\rtimes GL^+_d(\mathbb Q),\mathbb Z^d) $ is (isomorphic to) the full corner in the crossed
product 
$C_0(\mathcal A_f^d) \rtimes_\beta GL^+_d(\mathbb Q)$ corresponding to the projection 
$1_{\mathcal R^d} \in \mathcal A^d_f$. 
\end{proposition}
\begin{proof}
Observe first that the duality pairing of $\mathbb Q^d /\mathbb Z^d$ with $\mathcal R^d$ given by 
\[
\langle q , x\rangle = \exp (2\pi i \sum_{j=1}^d q_j x_j) 
\]
establishes a Gelfand-Fourier transform isomorphism of $C^*(\mathbb Q^d /\mathbb Z^d)$ to 
$C(\mathcal R^d)$. To see that the action $\beta$ restricted to $S$
is conjugate to the semigroup action of $S$ on $C^*(\mathbb Q^d / \mathbb Z^d)$
from part (ii) of \thmref{mainthm}, it is enough to carry out the computation on 
the characters $e(x) = \langle x,\cdot\rangle$. Thus 
\[
\alpha_s(e(x)) = R(s) \inv \sum_{[s^t y = x]} e(y)
\]
yields
\[
\hat\alpha_s(\langle x, \cdot\rangle) =  R(s) \inv \sum_{[s^t y = x]}\langle y, \cdot\rangle
=R(s) \inv \sum_{[s^t y = x]}\langle (s^t)\inv s^t y, \cdot\rangle
\]
 which, when  evaluated at $w\in \mathcal R^d$, equals $\langle  x, s\inv w\rangle$ or zero, according to
 $w\in s\mathcal R^d$ or not. Hence $\beta_s\circ i = i \circ \alpha_s$, where $i: C(\mathcal R^d)
\to C_0(\adel^d)$ is the obvious inclusion.

Since for every $s \in S$ there is a positive integer $m\in \mathbb N^\times$
such that $ms\inv \in S$, the semigroup $\mathbb N^\times$
of multiplicative positive integers is cofinal in $S$, and from this it follows
easily that $\cup_S \beta_{s\inv} C(\mathcal R^d)$ is dense in $C_0(\adel^d)$,
because $\cup_{m\in \mathbb N^\times} m\inv  \mathcal R^d $ is dense in $\adel^d$.
The result now follows from theorems 2.1 and 2.4 of \cite{minautex}.
\end{proof}
  \end{example}
\begin{example} 
Suppose $\K$ is a number field with ring of integers $\mathcal O$, and fix an integer $d \geq 1$.
Define
 \begin{align}
G &: = GL_d(\mathcal K),\notag\\
S &: = GL_d(\mathcal K) \cap \mathcal M_d(\mathcal O)\notag\\
\Gamma_0 &:= \mathcal O^d,\notag\\
N & : = \mathcal K^d \notag
\end{align}
and let $\Gamma := N \rtimes G$ where the action of $g\in G$ on $N$ is by $\psi_g(n) = (g^t)\inv n$.

\begin{proposition} 
With the above notation, $(\Gamma, \Gamma_0)$ is a Hecke pair, whose Hecke C*-algebra is 
isomorphic to the semigroup crossed product $C^*(\mathcal K^d / \mathcal O^d) \rtimes_\alpha S$,
and to a corner, corresponding to the 
characteristic function of the maximal compact subring $\mathcal R_\mathcal K^d$ 
of $\mathcal A_\mathcal K^d$, in the crossed product 
$C_0(\mathcal A_\mathcal K^d) \rtimes_\beta GL_d(\mathcal K)$,
where $\beta$ denotes the action of $GL_d(\mathcal K)$  on $C_0(\mathcal A_\mathcal K^d)$
given by $\beta_g (f) (x) = f(g\inv x)$.
\end{proposition}
The proof consists of an application of \thmref{mainthm} and of 
Theorems 2.1 and 2.4 of \cite{minautex}, and is entirely analogous to that of the preceding
example.
\end{example}

\begin{example}
 Motivated by Example 4.6 of \cite{bre} we consider the following situation.
 Let $R$ be a commutative unital ring and let $\Gamma_0$ be an $R$-module. 
Suppose $S \subset R$ is a multiplicatively closed set without zero-divisors,
 so that one has localizations  $S\inv R$ of the ring $R$ and $S\inv \Go$ of the module $\Go$
that contain, respectively, copies of $R$ and $\Go$.

Then $G = S\inv S$ is a subgroup of the unit group of $S\inv R$, there is an action
$\psi$ of $G$ on $N:= S\inv \Gamma_0$ given simply by $\psi_g(n) = g\inv n$, and we define $\Ga := N\rtimes G$.

\begin{proposition}
Assume $|s\inv \Gamma_0 / \Gamma_0| < \infty$ for each $s\in S$ and let $\Ga :=S\inv \Gamma_0\rtimes S\inv S$. 
Then  $(\Ga , \Go) $ is a Hecke pair and 
its Hecke algebra (respectively C*-algebra) is the semigroup crossed product
$\mathbb C (S\inv \Gamma_0 / \Gamma_0) \rtimes S)$ (respectively $C^*(S\inv \Gamma_0 / \Gamma_0) \rtimes S$).
\end{proposition}

\end{example}
\begin{example}
In the above examples $\Go$ is normal in $N$, so Corollary \ref{normal} is all that is needed. Next we 
borrow from \cite{lvf} an interesting example arising in number theory 
 in which the inclusion $\Go \subset N$ is not
normal, so that the full generality of theorems \ref{mainthm} and \ref{maincstarthm} is needed.
Let $\K$ be a number field, with ring of integers $\mathcal O$, and assume further that $\K$ is quadratic imaginary,
so that the unit group $\mathcal O^*$ is finite.
 Consider the group inclusion $\Go \subset \Ga$ with
 $\Go :=  \smatr{\mathcal O}{\mathcal O^*}$ and  $\Ga :=  \smatr{\K}{\K^*}$. The map 
$\smatr{y}{x} \mapsto  x\mathcal O^*$ is a
homomorphism from $\Ga$ onto the group $G := \K^*/\mathcal O^*$ 
of principal (fractional) ideals of $\K$; its  kernel is $N:=\smatr{\K}{\mathcal O^*}$. In general $\Go$ is 
not normal in $N$, but the double $\Go$--coset of $\smatr{y}{x} \in N$ is $\smatr{y\mathcal O^* +\mathcal
O}{\mathcal O^*}$, which contains at most  $|\mathcal O^*|$ right cosets, so $(N,\Go)$ is a Hecke pair.

 Let $S:= \mathcal O^\times /\mathcal O^*$ be the semigroup of principal {\em integral}
ideals, and let $\sigma: (\mathcal O^\times/\mathcal O^*) \to \mathcal O^\times$  be a cross section of
$\mathcal O^\times \to (\mathcal O^\times/\mathcal O^*)$ (see for example \cite[Lemma 2.5]{lvf}). Then  
$ (a \mathcal O^*) \mapsto \smatr{0}{\sigma(a\mathcal O^*)}$ determines a cross section $\tilde\sigma$ 
of the quotient  $\Ga \to \K^*/\mathcal O^*$, so we have that  $\Ga $ is a semidirect product 
\[
 \matr{\K}{\mathcal O^*} \rtimes (\K^*/\mathcal O^*),
\]
in which the action of $\K^*/\mathcal O^*$ is given by 
$\psi_g\smatr{y}{u} := \tilde\sigma(g) \smatr{y}{u} \tilde\sigma(g)\inv$.

For each  $ a\mathcal O^* \in S$ the subgroup 
$\tilde\sigma(a)\inv \Go \tilde\sigma(a) = \smatr{a\mathcal O}{\mathcal O^*}$ is of finite order in 
$\Go$ because $|\mathcal O / a\mathcal O| <
\infty$. It follows from \thmref{mainthm} that the pair
$ ({\K\rtimes\K^*}, {\mathcal O}\rtimes{\mathcal O^*} ) $
 is a Hecke pair and that 
\[
\mathcal H({\K\rtimes\K^*}, {\mathcal O}\rtimes{\mathcal O^*} ) \cong \mathcal H(\mathcal K \rtimes \mathcal O^*, \mathcal O
\rtimes \mathcal O^*)  \rtimes (\mathcal O^\times/\mathcal O^*).
\]
Moreover, the pair $ ({\K\rtimes\mathcal O^*}, {\mathcal O}\rtimes{\mathcal O^*} ) $ has a universal Hecke C*-algebra, and 
$ C^*_u(\Ga,\Go) \cong C^*_u(N,\Go) \rtimes (\mathcal O^\times/\mathcal O^*)$.
\begin{remark}
The assumption that $\mathcal K$ be quadratic imaginary is not necessary for the result,
but simplifies the argument;  we refer to \cite{lvf} for the proof of the general statement.
Here we limit ourselves to point out that $C^*_u(N,\Go)$ exists because its generators have finite spectrum and hence
uniformly bounded norm in any representation, and that 
the action of the principal integral ideals $\mathcal O^\times/\mathcal O^*$ by endomorphisms of 
$\mathcal H(\mathcal K \rtimes \mathcal O^*, \mathcal O \rtimes \mathcal O^*) $ is 
independent of the cross section $\sigma$, because the subgroup $\smatr{0}{\mathcal O^*}$ leaves $\Go$ invariant. 
\end{remark}
\end{example}


\begin{thebibliography}{20}

\bibitem{alr} J. Arledge, M. Laca and I. Raeburn,
{\em Semigroup crossed products and Hecke algebras arising from number
fields}, Documenta Math. {\bf 2} (1997), 115--138.

\bibitem{bin} M. W. Binder, {\em Induced factor representations of discrete
groups and their types}, J. Funct. Anal. {\bf 115} (1993), 294--312.

\bibitem{bla}  B. Blackadar, {\em Shape theory for $C^*$-algebras}, Math. Scand. 
{\bf 56} (1985),  249--275.

\bibitem{bos-con} J.-B. Bost and A. Connes, {\em Hecke algebras, Type
III factors  and phase transitions with spontaneous symmetry breaking
in number theory}, Selecta Math. (New Series) {\bf 1} (1995), 411--457.

\bibitem{bre} B. Brenken,
{\em Hecke algebras and semigroup crossed product $C^*$-algebras},
Pacific J. Math. {\bf 187} (1999), 241--262.

\bibitem{hal} R. W. Hall, {\it Hecke $C^*$-algebras}, Ph.D. thesis, The
  Pennsylvania State University, December 1999.

\bibitem{kri} A. Krieg, {Hecke Algebras}, Mem. Amer. Math. Soc. {\bf 87} (1990), No. 435.

\bibitem{quasilat} M. Laca and I. Raeburn, {\em Semigroup crossed products 
and the Toeplitz algebras of nonabelian groups}, 
J. Funct. Anal., {\bf 139} (1996), 415--440.

\bibitem{bcalg}  M. Laca and I. Raeburn, {\em  A semigroup crossed
product arising in number theory},  J. London Math. Soc., {\bf 59}
(1999), 330--344.

\bibitem{primbc} M. Laca and I. Raeburn, 
{\em The ideal structure of the Hecke $C^*$-algebra of Bost and Connes},
Math. Ann., {\bf 127} (2000), 433--451.


\bibitem{minautex} M. Laca, {\it From endomorphisms to automorphisms and  back:
 dilations and full corners}, J. London Math. Soc. (2) {\bf 61}
(2000), 893--904.

\bibitem{lvf}  M. Laca, M. van Frankenhuysen, {\em Phase transitions on Hecke C*-algebras and
class field theory}, in preparation.

\bibitem{LarsenR_MathScand} N. S. Larsen and I. Raeburn, {\em Faithful representations of
   crossed products by actions of $\mathbb{N}^k$}, Math. Scand., to appear.

\bibitem{lar-rae} N. S. Larsen and I. Raeburn, 
{\em Representations of Hecke algebras and dilations of semigroup crossed products},
preprint (June 2001).



\bibitem{mur}
G. J. Murphy, {\em Crossed products of $C^*$-algebras by endomorphisms},
Integral Equations \& Operator Theory {\bf 24} (1996), 298--319.


\bibitem{tza} K. Tzanev, {C*-alg\`{e}bres de Hecke et K-theorie}, Th\`ese de Doctorat, Universit\'{e} de Paris 7,
December 2000.

\end{thebibliography}
\end{document}